\documentclass[11pt]{elsart} \usepackage{amsmath,amssymb,latexsym}

\newenvironment{proof}[1][Proof.]{\begin{trivlist}
\item[\hskip \labelsep {\bfseries #1}]}{\end{trivlist}}

\newcommand{\openbox}{\leavevmode
  \hbox to.77778em{%
  \hfil\vrule
  \vbox to.675em{\hrule width.6em\vfil\hrule}%
  \vrule\hfil}}
\DeclareRobustCommand{\textsquare}{%
  \begingroup \usefont{U}{msa}{m}{n}\thr@@\endgroup
}

\def\QQfnmark#1{\footnotemark}

\newtheorem{lemma}{Lemma}

\newtheorem{corollary}{Corollary}

\newtheorem{theorem}{Theorem} 
\begin{document}

\title{Parameter-based Fisher's information of orthogonal polynomials}
\maketitle

\begin{frontmatter}
\author[Departamento1,Instituto]{J.\ S.\ Dehesa\corauthref{cor1}},\ead{dehesa@ugr.es}
\corauth[cor1]{Correspondence address: Instituto Carlos I de
F\'{\i}sica Te\'{o}rica y Computacional, Universidad de Granada,
E-18071 Granada, Spain}
\author[Departamento1,Instituto]{B. Olmos} \&
\author[Departamento2,Instituto]{R.J. Y\'{a}\~{n}ez}\ead{ryanez@ugr.es}

 \address[Departamento1]{Departamento de F\'{\i}sica Moderna,
   Universidad de Granada, 18071-Granada, Spain}

 \address[Departamento2]{Departamento de Matem\'atica Aplicada,
   Universidad de Granada, 18071-Granada, Spain}

 \address[Instituto]{Instituto Carlos I de F\'{\i}sica Te\'orica y
   Computacional, Universidad de Granada, 18071-Granada, Spain}

\begin{abstract}
  The Fisher information of the classical orthogonal polynomials with
  respect to a parameter is introduced, its interest justified and its
  explicit expression for the  Jacobi,
  Laguerre, Gegenbauer and Grosjean  polynomials found.
\end{abstract}
\end{frontmatter}
\newpage

\section{Introduction}

Let
$\{\rho_\theta(x)\equiv\rho(x|\theta);x\in\Omega\subset\mathbb{R}\}$
be a family of probability densities (or, more generally, likelihood
functions) parametrized by a parameter $\theta\in\mathbb{R}$. The
Fisher information of $\rho_\theta(x)$ with respect to the parameter
$\theta$ is defined \cite{CRA,FIS1,FIS2} as
\begin{eqnarray}\label{fisdef}
  I(\rho_\theta)&:=&\int_\Omega\left[\frac{\partial
      \ln\rho(x|\theta)}{\partial\theta}\right]^{2}\rho(x|\theta) 
dx=\int_\Omega \frac{\left[\frac{\partial
      \rho(x|\theta)}{\partial\theta} \right]^{2}}{\rho(x|\theta)}dx\nonumber \\
  &=&4\int_\Omega\left\{\frac{\partial \left[
        \rho(x|\theta)^{1/2}\right] }{\partial\theta}\right\}^{2}dx
\end{eqnarray}
This quantity refers to information about an unknown parameter in the
probability density calculated from observed atoms. It is a measure of
the ability to estimate the parameter $\theta$. It gives the minimum
error in estimating the parameter of the probability density
$\rho_\theta(x)$ \cite{COV}. The notion of Fisher information was
introduced by Sir R.A. Fisher in estimation theory \cite{FIS1,FIS2}.
Nowadays, it is being used in numerous scientific areas ranging from
statistics \cite{AMA,COX,FIS2,RAO,SCH}, information theory \cite{COV}
to signal analysis \cite{VIG} and quantum physics
\cite{FRI1,FRI2,FRI3,LUO2}. This information-theoretic quantity has,
among other characteristics, a number of important properties, beyond
the mere non-negativity, which deserve to be resembled here.

\begin{enumerate}

\item Additivity for independent events \cite{FIS2}. In the case that
  $\rho(x,y|\theta)=\rho_1(x|\theta)\cdot\rho_2(y|\theta)$,
  it happens that
  \[
  I\left[\rho(x,y|\theta)\right]=I\left[\rho_1(x|\theta)\right]\cdot
  I\left[\rho_2(y|\theta)\right]
  \]
\item Scaling invariance \cite{RAO}. The Fisher information is
  invariant under sufficient transformations $y=t(x)$, so that
  \[
  I\left[\rho(y|\theta)\right]=I\left[\rho(x|\theta)\right].
  \]
  This property is not only closely related to the Fisher maximum
  likelihood method \cite{RAO} but also it is very important for the
  theory of statistical inference.

\item Cramer-Rao inequality \cite{STA}. It states that the reciprocal
  of the Fisher information $I(\rho_\theta)$ bounds from below the
  mean square error of an unbiased estimator $f$ of the parameter
  $\theta$; i.e.
  \[
  \sigma^{2}(f)\ge\frac{1}{I(\rho_\theta)},
  \]
  where $\sigma^{2}(f)$ denotes the variance of $f$. This inequality,
  which lies at the heart of statistical estimation theory \cite{SCH},
  shows how much information the distribution provides about a
  parameter. Moreover, this says that the Fisher information
  $I(\rho_\theta)$ is a more sensitive indicator of the localization
  of the probability density than the Shannon entropy power \cite[p.
  240]{CAR}.

\item Relation to other information-theoretic properties. The Fisher
  information is related to the Shannon entropy of $\rho(x|\theta)$
  via the elegant de Bruijn's identity \cite{COV,LUO,STA}
  \[
  \frac{\partial}{\partial\theta}S(\tilde{\rho}_\theta)=\frac{1}{2}I(\tilde{\rho}_\theta),
  \]
  where $\tilde{\rho}_\theta$ denotes the convolution probability
  density of any probability density $\rho(x|\theta)$ with the
  normal density with zero mean and variance $\theta>0$, and
  $S(\tilde{\rho}_\theta):=-\int_\Omega
  \tilde{\rho}_\theta(x)ln\tilde{\rho}_\theta(x)dx$ is the Shannon
  entropy of $\tilde{\rho}_\theta(x)$.
 
  Moreover, the Fisher information $I(\rho_\theta)$ satisfies, under
  proper regularity conditions, the limiting property
  \cite{FRI1,FRI2,HAY}
  \[
  I(\rho_\theta)=\lim_{\epsilon \to
    0}\frac{2}{\epsilon^{2}}D\left(\rho_{\theta+\epsilon}\Vert\rho_\theta\right),
  \]
  where the symbol $D(p\Vert q):=\int_\Omega p(x)
  ln\frac{p(x)}{q(x)}dx$ denotes the relative entropy or
  Kullback-Leibler divergence of the probability densities $p(x)$ and
  $q(x)$. Further connections of the Fisher information with other
  information-theoretic properties are known; see e.g.
  \cite{BARR,DEM,FRI1,FRI2,LUT,MOH,SEA,STA}.

\item Applications in quantum physics. The Fisher information
  $I(\rho_\theta)$ plays a fundamental role in the quantum-mechanical
  description of physical systems
  \cite{CAR,DEH2,DEH3,FRI1,FRI2,FRI3,GON1,HAL1,HAL2,LUO,NAG,ROM1,ROM2,ROM4,RDY,ROM3,ROM5,SEA}.It
  has been shown

  \begin{enumerate}
  \item to be a measure of both disorder and uncertainty
    \cite{FRI1,FRI2,FRI3} as well as a measure of nonclassicality for
    quantum systems \cite{HAL1,HAL2},
  \item to describe, some factor apart, various macroscopic quantities
    such as, for example, the kinetic energy \cite{LUO,SEA,YAN} and
    the Weisz\"{a}cker energy \cite{ROM2,ROM4,RDY},
  \item to derive the Schr\"{o}dinger and Klein-Gordon equations of
    motion \cite{FRI1,FRI2} as well as the Euler equation of the
    density functional theory \cite{NAG}, from the principle of
    minimum Fisher information,
  \item to predict the most distinctive non-linear spectral phenomena,
    known as avoided crossings, encountered in atomic and molecular
    systems under strong external fields \cite{GON1}, and
  \item to be involved in numerous uncertainty inequalities
    \cite{DEH4,HAL2,LUO,LUO2,ROM4,ROM3,ROM5,STA}
  \end{enumerate}

\end{enumerate}

These applications are most apparent when $\theta$ is a parameter of
locality, so that $\rho_\theta(x)=\rho(x+\theta)$. Then, the Fisher
information for locality, also called intrinsic accuracy
\cite{COX,SEA}, does not depend on the parameter $\theta$; so, without
loss of generality, we may set the location at the origin and the
Fisher information of the density $\rho_0(x)\equiv\rho(x)$ becomes
\[
I(\rho)=\int\frac{\left[ \rho'(x)\right] ^2}{\rho(x)}dx,
\]
where $\rho'(x)$ denotes the first derivative of
$\rho(x)$. The locality Fisher information or Fisher information
associated with translations of an one-dimensional observable $x$ with
probability density $\rho(x)$ has been calculated in the literature
\cite{COV,VIG} for some simple families of probability densities. In
particular, it is well-known that the locality Fisher information of
the Gaussian density is equal to the reciprocal of its variance, what
illustrates the spreading character of these information-theoretic
measures. 

Recently, J.S. Dehesa and J. S\'{a}nchez-Ruiz \cite{SAN} have exactly
derived the locality Fisher information of a wider and much more
involved class of probability densities, the Rakhmanov densities,
defined by
\begin{equation}\label{Rakh}
  \rho_n(x)=\frac{1}{d_n^2}p_n^2(x)\omega(x)\chi_{\left[a,b\right]}(x),
\end{equation}
where $\chi_{\left[a,b\right]}(x)$ is the characteristic
function for the interval $\left[a,b\right]$, and $\{p_n(x)\}$ denotes
a sequence of real polynomials orthogonal with respect to the
nonnegative definite weight function $\omega(x)$ on the interval
$\left[a,b\right]\subseteq\mathbb{R}$, that is
\begin{equation}\label{dncuadrado}
  \int_a^bp_n(x)p_m(x)\omega(x)dx=d_n^2\delta_{n,m}
\end{equation}
with $\deg p_n(x)=n$. As first pointed out by E.A. Rakhmanov
\cite{RAK}, see also B. Simon \cite{SIM}, these densities play a
fundamental role in the analytic theory of orthogonal polynomials. In
particular, he has shown that these probabilty densities govern the
asymptotic behavior of the ratio $p_{n+1}(x)/p_n(x)$ as
$n\rightarrow\infty$. On the other hand, the Rakhmanov densities of
the classical orthogonal polynomials of a real continuous variable
describe the quantum-mechanical probability densities of ground and
excited states of numerous physical systems with an exactly solvable
Schr\"{o}dinger equation \cite{BAG,COO,NIK,USH}, particularly the most
common prototypes (harmonic oscillator, hydrogen atom,...), in
position and momentum spaces \cite{GAL}.

These two fundamental and applied reasons have motivated an increasing
interest for the determination of the spreading of the classical
orthogonal polynomials $\{p_n(x)\}$ throughout its interval of
orthogonality by means of the information-theoretic measures of their
corresponding Rakhmanov densities $\rho(x)$
\cite{BUY,DEH1,DEH2,DEH3,DEH4,DEH5,DVA,SAN,YAN}.

The Shannon information entropy of these densities has been examined
numerically \cite{BUY}. On the theoretical side, let us point out that
its asymptotics $(n\rightarrow\infty)$ is well known for all classical
orthogonal polynomials, but its exact value for every fixed $n$ is
only known for Chebyshev polynomials \cite{YVD} and some Gegenbauer
polynomials \cite{BUYA,SANL,SAR}. To this respect see \cite{DEH5}
which reviews the knowledge up to 2001. The variance and Fisher
information entropy of the Rakhmanov densities have been found in a
closed and compact form for all classical orthogonal polynomials
\cite{DEH3,SAN}. For other functionals of these Rakhmanov densities
see Ref. \cite{SAN1}.

In this paper we shall calculate the Fisher information of the real
and continuous classical orthogonal polynomials (Gegenbauer, Grosjean,
Jacobi, Laguerre) with respect to the parameter(s) of the polynomials.
We begin in Section 2 with the definition of this notion and, as well,
we collect here some basic properties of the classical orthogonal
polynomials which will be used later on. Then, the Fisher
information with respect to a parameter is fully determined for Jacobi
and Laguerre polynomials in Section 3, and for Gegenbauer and
Grosjean polynomials in Section 4. Finally, conclusions and some open
problems are given.

\section{Some properties of the parameter-dependent classical
  orthogonal polynomials}

Let $\{\tilde{y}_n(x;\theta)\}_{n\in \mathbb{N}_0}$ stand for the
sequence of polynomials orthonormal with respect to the nonnegative
definite weight function $\omega(x;\theta)$ on the real support
$(a,b)$, so that
\begin{equation}\label{orth}
  \int^b_a\tilde{y}_n(x;\theta)\tilde{y}_m(x;\theta)\omega(x;\theta)dx=\delta_{nm},
\end{equation}
with $\deg \tilde{y}_n=n$. Here we shall consider the celebrated
classical families of 
Hermite $H_n(x)$, Laguerre $L_n^{(\alpha)}(x)$ and Jacobi
$J_n^{(\alpha,\beta)}(x)$ polynomials. The normalized-to-unity density functions
$\tilde{\rho}_n(x;\theta)$ defined as
\begin{equation}\label{density}
  \tilde{\rho}_n(x;\theta)=\omega(x;\theta) \tilde{y}_n^2(x;\theta)
\end{equation}
are called for Rakhmanov densities \cite{RAK,SIM}.

Here we gather various properties of the parameter-dependent classical
orthogonal polynomials in a real and continuous variable (i.e.
Laguerre and Jacobi) in the form of two Lemmas, which shall be used
later on. The weight function of these polynomials can be written as
\begin{equation}\label{peso}
  \omega(x;\theta)=h(x)\left[ t(x)\right] ^{\theta}
\end{equation}
with
\begin{equation}\label{pesoL}
  h_L(x)=e^{-x}\quad and\quad t_L(x)=x
\end{equation}
for the Laguerre case, $L_n^{(\theta)}(x)$, and
\begin{equation}\label{pesoJ}
  h_J(x)=(1+x)^{\beta}\quad and\quad t_J(x)=1-x
\end{equation}
for the Jacobi case, $P_n^{(\theta,\beta)}(x)$.

%\subsection*{Lemma 1}
\begin{lemma}
The derivative of the orthonormal polynomial $\tilde{y}_n(x;\theta)$
with respect to the parameter $\theta$ is given by
\begin{equation}\label{lema1}
  \frac{\partial}{\partial{\theta}}\tilde{y}_n(x;\theta)=\sum^{n}_{k=0} 
\tilde{A}_k(\theta)\tilde{y}_k(x;\theta)  
\end{equation}
with
\begin{eqnarray}
  \tilde{A}_k(\theta)&=&\frac{d_k(\theta)}{d_n(\theta)}A_k(\theta)\mbox{\quad for\quad}
 k=0,1,\ldots,n-1 \label{Ak} \\ 
  \tilde{A}_n(\theta)&=&A_n(\theta)-\frac{1}{d_n(\theta)}\frac{\partial}{\partial\theta}
\left[ d_n(\theta)\right],\label{An}
\end{eqnarray}
where $d_m^2(\theta)$ denotes, according to Eq. (\ref{dncuadrado}),
the norm of the orthogonal polynomial $p_m(x)=y_m(x;\theta)$, and
$A_k(\theta)$ with $k=0,1,\ldots$ are the expansion coefficients of
the derivative of $y_m(x;\theta)$ in terms of the system
$\{y_m(x;\theta)\}$; i.e.,
\begin{equation}\label{partial}
  \frac{\partial}{\partial\theta}y_m(x;\theta)=\sum^{m}_{k=0}A_k(\theta)y_k(x;\theta).
\end{equation}
\end{lemma}

Both quantities $d_m(\theta)$ and $A_m(\theta)$ are known in the
literature for the Laguerre and Jacobi cases. Indeed, the norms for
the Laguerre $L_k^{(\alpha)}(x)$ and the Jacobi
$P_k^{(\alpha,\beta)}(x)$ polynomials \cite{NIK} are
\begin{eqnarray}
  \left[ d_k^{(L)}(\alpha)\right]^2&=&\frac{\Gamma(k+\alpha+1)}{k!},
  \label{dkl} \\
  \left[
    d_k^{(J)}(\alpha,\beta)\right]^2&=&\frac{2^{\alpha+\beta+1}\Gamma(k+\alpha+1)
\Gamma(k+\beta+1)}{k!(2k+\alpha+\beta+1)\Gamma(k+\alpha+\beta+1)},\label{dkj}
\end{eqnarray}
respectively.

On the other hand the expansion coefficients in Eq. (\ref{partial})
are known \cite{FRO,KOE1,KOE,RON} to have the form
\begin{equation}\label{AkL}
  A_k^{(L)}=\frac{1}{n-k}\mbox{\quad for\quad}k=0,1,\ldots,n-1\mbox{\quad and\quad}A_n^{(L)}=0
\end{equation}
for the Laguerre polynomials $L_n^{(\alpha)}$, and
\begin{eqnarray}\label{AkJa}
  A_k^{(J_\alpha)}&=&\frac{\alpha+\beta+1+2k}{(n-k)(\alpha+\beta+1+n+k)}\\\nonumber
  &&\qquad \times\frac{(\beta+k+1)_{n-k}}{(\alpha+\beta+k+1)_{n-k}};\quad k=0,1,\ldots,n-1,\\\label{AnJa}
  A_n^{(J_\alpha)}&=&\sum^{n-1}_{k=0}\frac{1}{\alpha+\beta+1+n+k}\\\nonumber
  &=&\psi(1+\alpha+\beta+2n)-\psi(1+\alpha+\beta+n),
\end{eqnarray}
for the Jacobi expansion of
$\frac{\partial}{\partial\alpha}P_n^{(\alpha,\beta)}(x)$, and
\begin{eqnarray}\label{AkJb}
  A_k^{(J_\beta)}&=&(-1)^{n-k}\frac{\alpha+\beta+1+2k}{(n-k)(\alpha+\beta+1+n+k)}\\\nonumber
  &&\qquad \times\frac{(\alpha+k+1)_{n-k}}{(\alpha+\beta+k+1)_{n-k}};\quad k=0,1,\ldots,n-1,\\\label{AnJb}
  A_n^{(J_\beta)}&=&\sum^{n-1}_{k=0}\frac{1}{\alpha+\beta+1+n+k}\\\nonumber
  &=&\psi(1+\alpha+\beta+2n)-\psi(1+\alpha+\beta+n),
\end{eqnarray}
for the Jacobi expansion of
$\frac{\partial}{\partial\beta}P_n^{(\alpha,\beta)}(x)$. The symbol
$\psi(x)=\frac{\Gamma'(x)}{\Gamma(x)}$ denotes the well-known digamma
function. It is worth to remark that the superindices $J_\alpha$ and
$J_\beta$ indicate the expansion coefficients for the derivatives of
the Jacobi polynomial with respect to the first and second parameters, respectively.

Then, taking into account Eqs. (\ref{Ak})-(\ref{An}) and Eqs.
(\ref{dkl}) and (\ref{AkL}), Lemma 1 says that the expansion
coefficients $\tilde{A}_k(\alpha)$ for the orthonormal Laguerre
polynomials $\tilde{L}_n^{(\alpha)}(x)$ are:
\begin{eqnarray}
  \tilde{A}_k^{(L)}&=&\frac{1}{n-k}\left[
    \frac{(k+1)_{n-k}}{(k+\alpha+1)_{n-k}}\right]^{1/2};\quad
  k=0,1,\ldots,n-1 \label{AkLtilde} \\
  \tilde{A}_n^{(L)}&=&-\frac{\psi(n+\alpha+1)}{2} \label{AnLtilde}
\end{eqnarray}
In a similar way the same Lemma together with Eqs.
(\ref{Ak})-(\ref{An}) and (\ref{dkj})-(\ref{AnJa}) has allowed us to
find the expressions
\begin{eqnarray}
  \tilde{A}_k^{(J_\alpha)}&=&\left[\frac{(k+\beta+1)_{n-k}
      (k+1)_{n-k}}{(k+\alpha+1)_{n-k}
      (k+\alpha+\beta+1)_{n-k}}\frac{2n+\alpha+\beta+1}{2k+\alpha+\beta+1}\right]^{1/2}\label{AkJatilde} 
    \\\nonumber 
  &&\times\frac{2k+\alpha+\beta+1}{(n-k)(n+k+\alpha+\beta+1)};\quad
  k=0,1,\ldots,n-1 \\
  \tilde{A}_n^{(J_\alpha)}&=&\frac{1}{2}\left[
    2\psi(2n+\alpha+\beta+1)-\psi(n+\alpha+\beta+1)\right.\label{AnJatilde} \\\nonumber 
  &&\left.-\psi(n+\alpha+1)-\ln{2}+\frac{1}{2n+\alpha+\beta+1}\right],
\end{eqnarray}
for the expansion coefficients $\tilde{A}_k^{(J_\alpha)}$ of the
derivative with respect to the parameter $\alpha$ of the Jacobi
polynomials $P_n^{(\alpha,\beta)}(x)$. Finally, the same procedure
with Eqs. (\ref{Ak})-(\ref{An}), (\ref{dkj}), (\ref{AkJb}) and
(\ref{AnJb}) leads to the values
\begin{eqnarray}
  \tilde{A}_k^{(J_\beta)}&=&\left[\frac{(k+\alpha+1)_{n-k}
      (k+1)_{n-k}}{(k+\beta+1)_{n-k}
      (k+\alpha+\beta+1)_{n-k}}\frac{2n+\alpha+\beta+1}{2k+\alpha+\beta+1}\right]^{1/2}\label{AkJbtilde} 
   \\\nonumber 
  &&\times(-1)^{n-k}\frac{2k+\alpha+\beta+1}{(n-k)(n+k+\alpha+\beta+1)};\quad k=0,1,\ldots,n-1 \\
  \tilde{A}_n^{(J_\beta)}&=&\frac{1}{2}\left[
    2\psi(2n+\alpha+\beta+1)-\psi(n+\alpha+\beta+1)\right.\label{AnJbtilde} \\\nonumber 
  &&\left.-\psi(n+\beta+1)-\ln{2}+\frac{1}{2n+\alpha+\beta+1}\right],
\end{eqnarray}
for the expansion coefficients $\tilde{A}_k^{(J_\beta)}$
($k=0,1,\ldots,n$) in Eq. (\ref{lema1}) of the $\beta$-derivative of
the Jacobi polynomials $P_n^{(\alpha,\beta)}(x)$.

%\subsection*{Lemma 2}
\begin{lemma}
The parameter-dependent classical orthonormal polynomials
$\tilde{y}_n(x;\theta)$ satisfy
\begin{eqnarray*}
  &(a)&\qquad\int_a^b\frac{\partial\omega(x;\theta)}{\partial\theta}\left[\tilde{y}_n(x;\theta)
  \right]^2dx=-2\tilde{A}_n(\theta)\\
  &(b)&\qquad\int_a^b\frac{\partial\omega(x;\theta)}{\partial\theta}\tilde{y}_n(x;\theta)
\tilde{y}_k(x;\theta)dx=-\tilde{A}_k(\theta)\quad;\quad k=0,1,\ldots,n-1\\ 
  &(c)&\qquad\int_a^b\frac{\partial^2\omega(x;\theta)}{\partial\theta^2}
\left[\tilde{y}_n(x;\theta)\right]^2dx=2\sum^{n}_{k=0}(\tilde{A}_k(\theta))^2
+2(\tilde{A}_n(\theta))^2-2\frac{\partial \tilde{A}_n(\theta)}{\partial\theta} 
\end{eqnarray*}
\end{lemma}
\begin{proof}
To prove the integrals (a) and (b) we have to derive with respect to
the parameter $\theta$ the orthonormalization condition (\ref{orth})
for $m=n$ and $m=k\neq n$, respectively. Then one has to use the Lemma
1 and again Eq. (\ref{orth}), and the results follow.
 
The integral (c) is obtained by deriving the integral (a) with respect
to $\theta$ and taking into account the values of the two previous
integrals (a) and (b). \hfill \openbox
\end{proof} 

\section{Parameter-based Fisher information of Jacobi and Laguerre
  polynomials}

The distribution of the orthonormal polynomials
$\tilde{y}_n(x;\theta)$ on their orthonormality interval and the
spreading of the associated Rakhmanov densities can be most
appropriately measured by means of their information-theoretic
measures, the Shannon entropy \cite{SHA} and the Fisher information
\cite{FIS1,FIS2}. The former has been theoretically
\cite{APT,DEH5,DVA} and numerically \cite{BUY} examined for general
orthogonal polynomials, while for the latter it has been studied the
Fisher information associated with translations of the variable (i.e.
the locality Fisher information) both analytically \cite{SAN} and
numerically \cite{DEH3}. Here we extend this study by means of the
computation of a more general concept, the parameter-based Fisher
information of the polynomials $\tilde{y}_n(x;\theta)$. This quantity
is defined as the Fisher information of the associated Rakhmanov
density (\ref{density}) with respect to the parameter $\theta$; that
is, according to Eq. (\ref{fisdef}), by
\begin{equation}\label{In}
  I_n(\theta)=4\int_a^b\left\lbrace\frac{\partial}{\partial\theta}
\left[ \tilde{\rho}_n(x;\theta)\right]^{1/2} \right\rbrace ^2dx 
\end{equation}
with
\[
\left[ \tilde{\rho}_n(x;\theta)\right]^{1/2}=\left[
  \omega(x;\theta)\right]^{1/2}\tilde{y}_n(x;\theta).
\]

%\subsection*{Theorem 1}
\begin{theorem}
The parameter-based Fisher information $I_n(\theta)$ of the
parameter-dependent classical orthonormal polynomials
$\tilde{y}_n(x;\theta)$ (i.e., Jacobi and Laguerre) defined by Eq.
(\ref{In}) has the value
\begin{equation}\label{Theo}
  I_n(\theta)=2\sum^{n-1}_{k=0}\left[
    \tilde{A}_k(\theta)\right]^2-2\frac{\partial\tilde{A}_n(\theta)}{\partial\theta}, 
\end{equation}
where $\tilde{A}_k(\theta)$, $k=0,1,\ldots,n$ are the expansion
coefficients of the derivative with respect to $\theta$ of
$\tilde{y}_n(x;\theta)$ in terms of the polynomials $\left\lbrace
  \tilde{y}_k(x;\theta)\right\rbrace_{k=0}^n$, which are given by
Lemma 1. See Eqs. (\ref{AkLtilde})-(\ref{AnLtilde}) and
(\ref{AkJatilde})-(\ref{AnJbtilde}) for the Laguerre and Jacobi
families, respectively.
\end{theorem}

\textbf{\underline{Remark.}} Let us underline that the Fisher
quantities of orthogonal, monic orthogonal and orthonormal polynomials
have the same value because of Eqs. (\ref{density}) and (\ref{In}),
keeping in mind the probabilistic character of the Rakhmanov density
and the fact that all these polynomials are orthogonal with respect to
the same weight function.

%\textbf{{Proof}}
\begin{proof}
To prove this theorem we begin with Eq. (\ref{In}). Then, the
derivative with respect to $\theta$ and Lemma 1 lead to
\begin{eqnarray*}
  \frac{\partial}{\partial\theta}\left[
    \tilde{\rho}_n(x;\theta)\right]^{1/2}&=&\left[
    \omega(x;\theta)\right]^{1/2}\sum^{n}_{k=0}\tilde{A}_k(\theta)
\tilde{y}_k(x;\theta)+\frac{\partial\left[
    \omega(x;\theta)\right]^{1/2}}{\partial\theta}
\tilde{y}_n(x;\theta). 
\end{eqnarray*}
The substitution of this expression into Eq. (\ref{In}) and the
consideration of the orthonormalization condition (\ref{orth}) have
led us to
\[
I_n(\theta)=J_1+J_2+J_3,
\]
where
\[
J_1=4\int_a^b\omega(x;\theta)\left(
  \sum^{n}_{k=0}\tilde{A}_k(\theta)\tilde{y}_k(x;\theta)\right)^2dx=4\sum^{n}_{k=0}\left[
  \tilde{A}_k(\theta)\right]^2,
\]
\begin{eqnarray*}
  J_2&=&4\int_a^b\left(\frac{\partial\left[ \omega(x;\theta)\right]
      ^{1/2}}{\partial\theta}\right)^2\left[
    \tilde{y}_n(x;\theta)\right]^2dx, 
\end{eqnarray*}
and
\begin{eqnarray*}
  J_3=8\sum^{n}_{k=0}\tilde{A}_k(\theta)\int_a^b\left[
    \omega(x;\theta)\right]^{1/2}\frac{\partial\left[
      \omega(x;\theta)\right]^{1/2}}{\partial\theta}\tilde{y}_n(x;\theta)\tilde{y}_k(x;\theta)dx. 
\end{eqnarray*}
Now we take into account that the weight function of the
parameter-dependent families of classical orthonormal polynomials in a
real and continuous variable (i.e. Laguerre and Jacobi) has the form
$\omega(x;\theta)=h(x)\left[t(x)\right]^{\theta}$, so that
\[
\left\lbrace \frac{\partial\left[ \omega(x;\theta)\right]
    ^{1/2}}{\partial\theta}\right\rbrace^2=\frac{1}{4}\omega(x;\theta)
\left[\ln{t(x)}\right]^2=\frac{1}{4}\frac{\partial^2\omega(x;\theta)}{\partial\theta^2} ,
\]
\[
\left[ \omega(x;\theta)\right]^{1/2}\frac{\partial\left[
    \omega(x;\theta)\right]^{1/2}}{\partial\theta}=\frac{1}{2}\omega(x;\theta)
\ln{t(x)}=\frac{1}{2}\frac{\partial\omega(x;\theta)}{\partial\theta}.
\]
The use of these two expressions in the integrals $J_2$ and $J_3$
together with Lemma 2 leads to Eq. (\ref{Theo}). \hfill \openbox
\end{proof}
%\begin{flushright}
%  $\blacksquare$
%\end{flushright}

%\subsection*{Corollary 1}
\begin{corollary}
The Fisher information with respect to the parameter $\alpha$,
$I_n^{(L)}(\alpha)$, of the Laguerre polynomial
$\tilde{L}_n^{(\alpha)}(x)$ is
\begin{eqnarray}
  I_n^{(L)}(\alpha)&=&2\sum^{n-1}_{k=0}\left[
    \tilde{A}_k^{(L)}\right]^2-2\frac{\partial
    \tilde{A}_n^{(L)}}{\partial\alpha}\\\nonumber 
  &=&\psi^{(1)}(n+\alpha+1)+\frac{2n}{n+\alpha}{}_4F_3\left(\begin{array}{c|c}
      1\quad1\quad1\quad1-n & {}\\
      {} & 1 \\
      2\quad2\quad1-\alpha-n & {} \\
    \end{array}\right),
\end{eqnarray}
where $\psi^{(1)}(x)=\frac{d}{dx}\psi(x)$ is the trigamma function.
\end{corollary}
%\subsection*{Corollary 2}
\begin{corollary}
The Fisher information with respect to the parameter $\alpha$,
$I_n^{(J_\alpha)}(\alpha,\beta)$, of the Jacobi polynomial
$\tilde{P}_n^{(\alpha,\beta)}(x)$ has the value
\begin{eqnarray}
  I_n^{(J_\alpha)}(\alpha,\beta)&=&2\sum^{n-1}_{k=0}\left[
    \tilde{A}_k^{(J_\alpha)}\right]^2-2\frac{\partial
    \tilde{A}_n^{(J_\alpha)}}{\partial\alpha}\\\nonumber 
  &=&2\frac{\Gamma(n+\beta+1)n!(2n+\alpha+\beta+1)}{\Gamma(n+\alpha+1)
\Gamma(n+\alpha+\beta+1)}\\\nonumber 
  &&\times\sum^{n-1}_{k=0}\frac{\Gamma(k+\alpha+1)\Gamma(k+\alpha+\beta+1)(2k+\alpha+\beta+1)}
{\Gamma(k+\beta+1)k!(n-k)^2(n+k+\alpha+\beta+1)^2}\\\nonumber 
  &&-2\psi^{(1)}(2n+\alpha+\beta+1)+\psi^{(1)}(n+\alpha+\beta+1)\\\nonumber
  &&+\psi^{(1)}(n+\alpha+1)+\frac{1}{(2n+\alpha+\beta+1)^2},
\end{eqnarray}
and the Fisher information with respect to the parameter $\beta$,
$I_n^{(J_\beta)}(\alpha,\beta)$, of the Jacobi polynomial
$\tilde{P}_n^{(\alpha,\beta)}(x)$ has the value
\begin{eqnarray}
  I_n^{(J_\beta)}(\alpha,\beta)&=&2\sum^{n-1}_{k=0}\left[
    \tilde{A}_k^{(J_\beta)}\right]^2-2\frac{\partial
    \tilde{A}_n^{(J_\beta)}}{\partial\beta}\\\nonumber 
  &=&2\frac{\Gamma(n+\alpha+1)n!(2n+\alpha+\beta+1)}{\Gamma(n+\beta+1)
\Gamma(n+\alpha+\beta+1)}\\\nonumber
  &&\times\sum^{n-1}_{k=0}\frac{\Gamma(k+\beta+1)\Gamma(k+\alpha+\beta+1)(2k+\alpha+\beta+1)}
{\Gamma(k+\alpha+1)k!(n-k)^2(n+k+\alpha+\beta+1)^2}\\\nonumber
  &&-2\psi^{(1)}(2n+\alpha+\beta+1)+\psi^{(1)}(n+\alpha+\beta+1)\\\nonumber
  &&+\psi^{(1)}(n+\beta+1)+\frac{1}{(2n+\alpha+\beta+1)^2}.
\end{eqnarray}
\end{corollary}
Both corollaries follow from Theorem 1 in a straightforward manner by
taking into account the expressions (\ref{AkLtilde})-(\ref{AnJbtilde})
for the expansion coefficients $\tilde{A}_k$ of the corresponding
families. Let us underline that $J_\alpha$ and $J_\beta$ indicate
Fisher informations with respect to the first and second parameter,
respectively, of the Jacobi polynomial.

\section{Parameter-based Fisher information of the Gegenbauer and
  Grosjean polynomials}

In this section we describe the Fisher information of two important
subfamilies of the Jacobi polynomials $P_n^{(\alpha,\beta)}(x)$: the
ultraspherical or Gegenbauer polynomials \cite{ASS,DEH5,NIK}, which
have $\alpha=\beta$, and the Grosjean polynomials of the first and
second kind \cite{DET,GRO,RONV,RDZ}, which have $\alpha+\beta=\pm1$,
respectively. Let us remark that the parameter-based Fisher information
for these subfamilies cannot be obtained from the expressions of the
similar quantity for the Jacobi polynomials (given by corollary 2) by
means of a mere substitution of the parameters, because it depends on
the derivative with respect to the  parameter(s) and now $\alpha$ and
$\beta$ are correlated.

The Gegenbauer polynomials $C_n^{(\lambda)}(x)$ are Jacobi-like
polynomials satisfying the orthogonality condition (\ref{dncuadrado})
with the weight function
$\omega_C(x;\lambda)=(1-x^2)^{\lambda-\frac{1}{2}}$, $\lambda >
-\frac{1}{2}$, and the normalization constant
\[
\left[d_k^{(C)}(\lambda)\right]^2=\frac{\pi2^{1-2\lambda}\Gamma(k+2\lambda)}{k!(k+\lambda)
\Gamma^2(\lambda)}
\]
so that they can be expressed as
\[
C_n^{(\lambda)}(x)=\frac{(2\lambda)_n}{(\lambda+\frac{1}{2})_n}P_n^{(\lambda-\frac{1}{2},
\lambda-\frac{1}{2})}(x)
\]
It is known \cite{KOE1,KOE} that the expansion (\ref{partial}) for the
derivative of $C_n^{(\lambda)}(x)$ with respect to the parameter
$\lambda$ has the coefficients
\begin{eqnarray*}
  A_k^{(C)}(\lambda)&=&\frac{2\left(1+(-1)^{n-k}\right)(k+\lambda)}{(k+n+2\lambda)(n-k)}
\quad\mbox{for}\quad k=0,1,\ldots,n-1\\ 
  A_n^{(C)}(\lambda)&=&\sum^{n-1}_{k=0}\frac{2(k+1)}{(2k+2\lambda+1)(k+2\lambda)}
+\frac{2}{k+n+2\lambda}\\
  &=&\psi(n+\lambda)-\psi(\lambda)
\end{eqnarray*}
Then, according to Eqs. (\ref{Ak}-\ref{An}) of Lemma 1, the expansion
(\ref{lema1}) for the derivative of the orthonormal Gegenbauer
polynomials has the following coefficients
\begin{eqnarray*}
  \tilde{A}_k^{(C)}(\lambda)&=&\left[\frac{\Gamma(k+2\lambda)n!(n+\lambda)}{\Gamma(n+2\lambda)
k!(k+\lambda)}\right]^{1/2}\\
  &&\times\frac{2\left(1+(-1)^{n-k}\right)(k+\lambda)}{(k+n+2\lambda)(n-k)}\quad\mbox{for}
\quad k=0,1,\ldots,n-1\\
  \tilde{A}_n^{(C)}(\lambda)&=&\psi(n+\lambda)-\psi(n+2\lambda)+\ln{2}+\frac{1}{2(n+\lambda)}.
\end{eqnarray*}
Theorem 1 provides, according to Eq. (\ref{Theo}), the following value
for the Fisher information of the Gegenbauer polynomials
$C_n^{(\lambda)}(x)$ with respect to the parameter $\lambda$:
\begin{eqnarray*}
  I_n^{(C)}(\lambda)&=&\frac{16n!(n+\lambda)}{\Gamma(n+2\lambda)}\sum^{n-1}_{k=0}
\frac{\left(1+(-1)^{n-k}\right)\Gamma(k+2\lambda)(k+\lambda)}{k!(k+n+2\lambda)^2(n-k)^2}\\
  &&-2\psi^{(1)}(n+\lambda)+4\psi^{(1)}(n+2\lambda)+\frac{1}{(n+\lambda)^2}.
\end{eqnarray*}
Let us now do the same job for the Grosjean polynomials of the first
and second kind, which are the monic Jacobi polynomials
$\hat{P}_n^{(\alpha,\beta)}(x)$ with $\alpha+\beta=\mp1$,
respectively. So, we have \cite{GRO,RONV}
\[
G_n^{(\alpha)}(x)=c_nP_n^{(\alpha,-1-\alpha)}(x)\quad,\quad-1<\alpha<0,
\]
and
\[
g_n^{(\alpha)}(x)=e_nP_n^{(\alpha,1-\alpha)}(x)\quad,\quad-1<\alpha<2,
\]
for the Grosjean polynomials of first and second kind, respectively,
with the values
\[
c_n=2^n\binom{2n-1}{n}^{-1}\quad,\quad e_n=2^n\binom{2n+1}{n}^{-1}.
\]
The Grosjean polynomials of the first kind $G_n^{(\alpha)}(x)$ satisfy
the orthogonality condition (\ref{dncuadrado}) with the weight
function
\[
\omega_G(x;\alpha)=\left(\frac{1-x}{1+x}\right)^\alpha\frac{1}{1+x},
\]
and the normalization constant
\[
\left[d_n^{(G)}(\alpha)\right]^2=\frac{2^{2n-1}\Gamma^2(n)}{\Gamma^2(2n)}
\Gamma(n+\alpha+1)\Gamma(n-\alpha).
\]
These polynomials, together with the Chebyshev polynomials of the
first, second, third and fourth kind, are the only Jacobi polynomials
for which the associated polynomials are again Jacobi polynomials \cite{GRO}.

Now, the expansion (\ref{partial}) for the derivative of the
polynomials $G_n^{(\alpha)}(x)$ with respect to the parameter $\alpha$
can be obtained as
\begin{eqnarray*}
  \frac{\partial G_n^{(\alpha)}(x)}{\partial\alpha}&=&\frac{\partial 
\hat{P}_n^{(\alpha,-1-\alpha)}(x)}{\partial\alpha}=\frac{\partial
\hat{P}_n^{(\alpha,\beta)}(x)}{\partial\alpha}\Bigg|_{\beta=-1-\alpha}-\frac{\partial
\hat{P}_n^{(\alpha,\beta)}(x)}{\partial\beta}\Bigg|_{\beta=-1-\alpha}
\\ 
  &=&\sum^{n-1}_{k=0}A_k^{(G)}(\alpha)G_n^{(\alpha)}(x)
\end{eqnarray*}
with
\[
A_k^{(G)}(\alpha)=\frac{2^{n-k+1}k}{n^2-k^2}\frac{\Gamma(2k)\Gamma(n+1)}
{\Gamma(2n)\Gamma(k+1)}\left[(k-\alpha)_{n-k}-(-1)^{n-k}(k+\alpha+1)_{n-k}\right]
\]
for $k=0,1,\ldots,n-1$ and $A_n^{(G)}(\alpha)=0$. Then, Lemma 1
provides the expansion (\ref{lema1}) for the derivative of the
orthonormal Grosjean polynomials with the coefficients
\begin{eqnarray*}
  \tilde{A}_k^{(G)}(\alpha)&=&\frac{2n}{n^2-k^2}\frac{(k-\alpha)_{n-k}-(-1)^{n-k}
(k+\alpha+1)_{n-k}}{\left[(k+\alpha+1)_{n-k}(k-\alpha)_{n-k}\right]^{1/2}}\\
  &&\qquad \qquad \qquad \mbox{for}\quad k=0,1,\ldots,n-1,\\
  \tilde{A}_n^{(G)}(\alpha)&=&\frac{1}{2}\left[\psi(n-\alpha)-\psi(n+\alpha+1)\right].
\end{eqnarray*}
Finally, Eq. (\ref{Theo}) of Theorem 1 allows us to find the following
value for the Fisher information of the Grosjean polynomials of the
first kind:
\begin{eqnarray*}
  I_n^{(G)}(\alpha)&=&8n^2\sum^{n-1}_{k=0}\frac{1}{\left(n^2-k^2\right)^2}
\frac{\left[(k-\alpha)_{n-k}-(-1)^{n-k}(k+\alpha+1)_{n-k}\right]^2}{(k+\alpha+1)_{n-k}(k-\alpha)_{n-k}}\\ 
  &&+\psi^{(1)}(n-\alpha)+\psi^{(1)}(n+\alpha+1).
\end{eqnarray*}
On the other hand, the Grosjean polynomials of the second kind
$g_n^{(\alpha)}(x)$ satisfy the orthogonality property
(\ref{dncuadrado}) with the weight function
\[
\omega_g(x;\alpha)=\left(\frac{1-x}{1+x}\right)^\alpha(1+x),
\]
and the normalization constant
\[
\left[d_n^{(g)}(\alpha)\right]^2=\frac{2^{2n+1}\Gamma^2(n)}{\Gamma^2(2n+2)}
\Gamma(n+\alpha+1)\Gamma(n-\alpha+2).
\]
Moreover, the derivative of these polynomials with respect to the
parameter $\alpha$ can be expanded in the form (\ref{lema1}) as
\begin{eqnarray*}
  \frac{\partial g_n^{(\alpha)}(x)}{\partial\alpha}&=&\frac{\partial
    \hat{P}_n^{(\alpha,1-\alpha)}(x)}{\partial\alpha}=\frac{\partial
    \hat{P}_n^{(\alpha,\beta)}(x)}{\partial\alpha}\Bigg|_{\beta=1-\alpha}
-\frac{\partial \hat{P}_n^{(\alpha,\beta)}(x)}{\partial\beta}\Bigg|_{\beta=1-\alpha}\\ 
  &=&\sum^{n-1}_{k=0}A_k^{(g)}(\alpha)g_n^{(\alpha)}(x)
\end{eqnarray*}
with
\begin{eqnarray*}
  A_k^{(g)}(\alpha)&=&\frac{2^{n-k+1}(k+1)}{(n-k)(n+k+2)}\frac{\Gamma(2k+2)\Gamma(n+1)}
{\Gamma(2n+2)\Gamma(k+1)}\\
  &&\times\left[(k+2-\alpha)_{n-k}-(-1)^{n-k}(k+\alpha+1)_{n-k}\right]
\end{eqnarray*}
for $k=0,1,\ldots,n-1$ and $A_n^{(g)}(\alpha)=0$. Then, Lemma 1 is
able to provide the analogous expansion (\ref{lema1}) for the
orthonormal polynomials with the coefficients
\begin{eqnarray*}
  \tilde{A}_k^{(g)}(\alpha)&=&\frac{2(k+1)}{(n-k)(n+k+2)}\frac{(k+2-\alpha)_{n-k}
-(-1)^{n-k}(k+\alpha+1)_{n-k}}{\left[(k+\alpha+1)_{n-k}(k+2-\alpha)_{n-k}\right]^{1/2}};\\
  &&\qquad \qquad \qquad \mbox{for}\quad k=0,1,\ldots,n-1,\\
  \tilde{A}_n^{(g)}(\alpha)&=&\frac{1}{2}\left[\psi(n+2-\alpha)-\psi(n+\alpha+1)\right].
\end{eqnarray*}
Finally, we obtain by means of Eq. (\ref{Theo}) of Theorem 1 the
Fisher information of the Grosjean polynomials of the second kind,
which turns out to have the value
\begin{eqnarray*}
  I_n^{(g)}(\alpha)&=&8\sum^{n-1}_{k=0}\frac{(k+1)^2}{(n-k)^2(n+k+2)^2}\\
  &&\times\frac{\left[(k+2-\alpha)_{n-k}-(-1)^{n-k}(k+\alpha+1)_{n-k}\right]^2}
{(k+\alpha+1)_{n-k}(k+2-\alpha)_{n-k}}\\
  &&+\psi^{(1)}(n+2-\alpha)+\psi^{(1)}(n+\alpha+1).
\end{eqnarray*}

\section{Conclusions and open problems}

In summary, we have calculated the parameter-based Fisher information
for the classical orthogonal polynomials of a continuous and real
variable with a parameter dependence; namely, the Jacobi and Laguerre
polynomials. Then we have evaluated the corresponding Fisher quantity
for the two most relevant parameter-dependent Jacobi polynomials
$P_n^{(\alpha,\beta)}(x)$: the Gegenbauer ($\alpha=\beta$) and the
Grosjean ($\alpha+\beta=\pm1$) polynomials.

This paper, together with Ref. \cite{SAN}, opens the way for the
developement of the Fisher estimation theory of the Rakhmanov density
for continuous and discrete orthogonal polynomials in and beyond the
Askey scheme. This fundamental task in approximation theory includes
the determination of the spreading of the orthogonal polynomials
throughout its orthogonality domain by means of the Fisher information
with a locality property. All these mathematical questions have a
straightforward application to quantum systems because their
wavefunctions are often controlled by orthogonal polynomials, so that
the probability densities which describe the quantum-mechanical states
of these physical systems are just the Rakhmanov densities of the
corresponding orthogonal polynomials.

Among the open problems let us first mention the computation of the
parameter-based Fisher information of the generalized Hermite
polynomials, the Bessel polynomials and the Pollaczek polynomials
\cite{CHI}. A much more ambitious problem is the evaluation of the
Fisher quantity for the general Wilson orthogonal polynomials
\cite{GEA}.
On the other hand, nothing is known for discrete orthogonal polynomials. In this
case, however, the very notion of the parameter-based Fisher
information is a subtle question.

\section{Acknowledgements}

This work has been partially supported by the MEC Project No.
FIS2005-00973, by the European Research Network on Constructive
Approximation (NeCCA) INTAS-03-51-6637 and by the J.A. Project of
Excellence with ref. FQM-481. We belong to the P.A.I. Group FQM-207 of
the Junta de Andaluc\'{i}a, Spain.

\bibliographystyle{ppcf} \bibliography{general}

\end{document}